\documentclass{birkjour}
\usepackage{amssymb}
\newcommand{\st}{\substack}
\newtheorem{thm}{Theorem}[section]
\newtheorem{cor}[thm]{Corollary}

\theoremstyle{definition}

\theoremstyle{remark}
\newtheorem{rem}[thm]{Remark}

\numberwithin{equation}{section}

\def\be{\begin{equation}}
\def\ee{\end{equation}}
\newcommand{\bthm}{\begin{thm}}
\newcommand{\ethm}{\end{thm}}
\newcommand{\bcor}{\begin{cor}}
\newcommand{\ecor}{\end{cor}}

\newcommand{\beq}{\begin{eqnarray}}
\newcommand{\beqq}{\begin{eqnarray*}}
\newcommand{\eeq}{\end{eqnarray}}
\newcommand{\eeqq}{\end{eqnarray*}}
\newcommand{\ba}{\begin{array}}
\newcommand{\ea}{\end{array}}
\newcommand{\bee}{\begin{enumerate}}
\newcommand{\eee}{\end{enumerate}}
\newcommand{\bpf}{\begin{proof}}
\newcommand{\epf}{\end{proof}}

\newcommand{\IN}{{\mathbb N}}

\newcounter{alphabet}
\newcounter{tmp}

\begin{document}

\title[Inequalities for Mertens functions]{Inequalities for weighted sums of Mertens functions}

\author[R. Balasubramanian]{Ramachandran Balasubramanian}
\address{Institute of Mathematical Sciences,\\
IV Cross Road, CIT Campus, \\
Taramani, Chennai 600 113, India.}
\email{balu@imsc.res.in}


\author[S. Ponnusamy]{Saminathan Ponnusamy
}
\address{Department of Mathematics,\\
Indian Institute of Technology Madras,\\
 Chennai-600 036, India.}
\email{samy@iitm.ac.in}

\author[K.-J. Wirths]{Karl-Joachim Wirths}
\address{Institut f\"ur Analysis und Algebra,\\
 TU Braunschweig, \\
38106 Braunschweig, Germany.}
\email{kjwirths@tu-bs.de}

\date{February 08, 2019}

\begin{abstract}
In this article we derive some polynomial inequalities for Mertens functions.
\end{abstract}
\maketitle

\section{Introduction and Main results}
Let $\mu(n)$ be the M\"obius function of the positive integer $n$, that is,
\bee
\item[{\rm (a)}] $\mu(1)\,=\,1$,
\item[{\rm (b)}] $\mu(n)\,=\,0$, if a square number is a divisor of $n$,
\item[{\rm (c)}] $\mu(n)\,=\,(-1)^r$, if $n$ is the product of $r$ pairwise disjoint prime numbers.
\eee
Suppose further that
\[
M(N)\,=\,\sum_{n=1}^N\mu(n)
\]
denotes the Mertens function.

During their efforts to prove a coefficient conjecture (see \cite[Conjecture 1]{OPW}) for some classes of univalent functions,
the second and the third authors of the present paper considered an inequality that concerned the Mertens function.
See also \cite{PW}. Notwithstanding the fact that these efforts had no success hitherto, the authors think that this
inequality and its proof are of independent interest and we want to present them in the sequel.

\bthm\label{th1}
Let $n\in \IN \setminus\{1\}$,  and as usual
\[\left[\frac{n}{k}\right]\,=\,\max\left\{j: j\in \IN, j\,\leq \,\frac{n}{k}\right\}.
\]
For $\lambda \in [0,1]$  the inequality
\begin{equation}\label{f12}
\sum_{j=0}^{n}\lambda^j\,-\,\sum_{k=1}^{n-1}\left\{M\left(\left[\frac{n-1}{k}\right]\right)
\left(\sum_{j=0}^{n-k}\lambda^j\right)\right\}\,=:\,\sum_{j=0}^nd_j\lambda^j\,\geq \,0
\end {equation}
 is valid. Equality occurs if and only if $\lambda = 0.$
\ethm
\bpf
For $\lambda =0$, we have the well known equation
\begin{equation}\label{f8}
1\,=\,\sum_{k=1}^nM\left(\left[\frac{n}{k}\right]\right),
\end{equation}
which shows that the assertion is valid in this case. We do not know the eldest reference for
(\ref{f8}), but we found that it has been proved and used to compute values of $M$ in
\cite{L,N}. Hence, we have to prove that (\ref{f12}) is valid with strict inequality $>$ instead of $\geq$ for $\lambda \in (0,1].$

Next, we consider the cases $2\,\leq n \,\leq 94$. Let us use the abbreviations $m=n-1$ and
\[ B_k\,=\,\sum_{j=0}^{m+1-k}\lambda^j ~\mbox{ for }~0\,\leq k \,\leq m+1.
\]
 As $\left[\frac{m}{k}\right]\,=\,1$ for $\left[\frac{m}{2}\right]+1\,\leq\,k\,\leq\,m,$ we have
$$
M\left(\left[\frac{m}{k}\right]\right)\,=\,\mu(1)\,=\,1~\mbox{ for }~\left[\frac{m}{2}\right]+1\,\leq\,k\,\leq\,m.
$$
Further it is known that
$$
M(j)\,\leq \,0 ~\mbox{ for }~ 2\,\leq\,j\,\leq 93.
$$
Since $B_k, 0\leq k \leq\,m+1,$ is monotonically decreasing, we get from the above and (\ref{f8}) the validity of (\ref{f12})
in the following way:
\beq
B_0\,-\,\sum_{k=1}^mM\left(\left[\frac{m}{k}\right]\right)B_k &= & B_0\,+\sum_{k=1}^{\left[\frac{m}{2}\right]}\left|M\left(\left[\frac{m}{k}\right]\right)\right|B_k\,-\,\sum_{k=\left[\frac{m}{2}\right]+1}^mB_k\nonumber\\
&\geq& B_0\,-\,\left(\sum_{k=1}^mM\left(\left[\frac{m}{k}\right]\right)\right)B_{\left[\frac{m}{2}\right]+1}\nonumber\\
&= &B_0\,-\,\,B_{\left[\frac{m}{2}\right]+1}\,> \,0. \nonumber
\eeq

It remains to consider the cases $\lambda \in (0,1]$ and $n\,\geq \,95$.

From now on we will use the abbreviation $M\left(\left[\frac{n-1}{r}\right]\right)\,=\,f(r)$ to make the formulas more readable
and now and then we use the abbreviation $m\,=\,n-1$. It is immediately seen that the coefficients $d_j$, $j=0, \ldots ,n,$ can be
calculated as follows:
\[ d_0\,=\,0,\quad d_n\,=\,1,
\]
and
\begin{equation}\label{f13}
d_{n-j}\,=\,1\,-\,\sum_{r=1}^jf(r)\,=\,\sum_{r=j+1}^{n-1}f(r) ~\mbox{ for }~
1\,\leq \,j\,\leq n-1,\end{equation}
as follows from (\ref{f8}). Formula (\ref{f13}) is equivalent to
\begin{equation}\label{f14}
d_j\,=\,\sum_{r=n-j+1}^{n-1}f(r)~\mbox{ for }~ 1\,\leq \,j\,\leq n-1.
\end{equation}
We begin our discussion by proving two items for coefficients $d_j$ with ``small'' indices. Firstly, we derive from (\ref{f14}) that
\begin{equation}\label{f15}
d_j\,=\,j-1 ~\mbox{ for }~ 2\,\leq\,j\,\leq \,n\,-\, \left[\frac{n-1}{2}\right].
\end{equation}
Since
\[ \frac{n}{2}\,+\,1\,\geq \,n\,-\, \left[\frac{n-1}{2}\right]\,\geq \frac{n+1}{2}, 
\]
we get
$$
\sum_{j=1}^{\frac{n}{2}}d_j\,\geq\, \frac{1}{2}\left(n\,-\, \left[\frac{n-1}{2}\right]\,-\,1\right)\left(n\,-\, \left[\frac{n-1}{2}\right]\,-\,2\right)
\,\geq \,\frac{n^2-4n+3}{8}.
$$
Further, we show that the inequalities
\[
d_j\,\geq\,0~\mbox{ for }~ 0\,\leq \,j\,\leq n\,-\,\left[\frac{n-1}{10}\right]
\]
are valid. Especially, we will use that this implies
\[
d_j\geq\,0~\mbox{ for }~ 0\,\leq \,j\,\leq \,\frac{9n}{10}.
\]
According to (\ref{f15}), we have to take into account the indices
\[
j\,\geq \, n\,-\,\left[\frac{n-1}{2}\right]\,+\,1.\] Since $M(2)\,=\,0$,   from (\ref{f14}) and (\ref{f15}), we get that
\[
d_j\,=\, n\,-\,1\,-\,\left[\frac{n-1}{2}\right]\,\geq \,\frac{n-1}{2}~\mbox{ for }~n\,-\,\left[\frac{n-1}{2}\right]\,+\,1\,\leq\,j
\,\leq\,n\,-\,\left[\frac{n-1}{3}\right].
\]
For
\[ n\,-\,\left[\frac{n-1}{3}\right]\,+\,1\,\leq\,j\,\leq\,n\,-\,\left[\frac{n-1}{5}\right]\]
we use $M(3)\,=\,M(4)\,=\,-1$ to achieve
\beqq
d_j &\geq &d_{n-\left[\frac{n-1}{5}\right]}\\
&\geq&\frac{n-1}{2}\,-\,\left(\left[\frac{n-1}{3}\right]\,-\,\left[\frac{n-1}{5}\right]\right)\\
&\geq& \frac{n-1}{2}\,-\,\frac{n-1}{3}\,+\,\frac{n}{5}\,-\,1\\
&=&\frac{11n\,-\,35}{30} \,\geq \,0,\,\, \mbox{ if $\,n\geq 4.$}
\eeqq
It is not difficult to continue in this way. At the end one arrives at the inequality
\[
d_j\,\geq \,d_{n-\left[\frac{n-1}{10}\right]}\,\geq\,\frac{40n\,-\,1116}{210}\,\geq \,0\]
for
\[
n\,\geq\,28,\,\,\textrm{and}\,\, n-\left[\frac{n-1}{2}\right]\,\leq \,j\,\leq\,n-\left[\frac{n-1}{10}\right].\]
Since $d_n\,=\,1$, we may restrict ourselves to prove
\[
S(\lambda)\,=\,\sum_{j=1}^md_j\lambda^j\,> \,0.
\]
To this end, let
\beqq
A&=&\left \{j: j \leq \frac{n}{2}\right \},\\
B&=&\left \{j:\frac{n}{2}\,<\, j \leq \frac{9n}{10}\right \},\\
C&=&\{j:0.9n \,<\, j \leq m, d_j \geq 0\}, ~\mbox{ and }\\
D&=&\{j:0.9n \,< j \leq m, d_j<0\}.
\eeqq
Then $A \cup B \cup C \cup D$ gives a partition of the set $\{j:1 \leq j \leq m\}$.

Let $a=\sum_{j \in A} d_j$,  $ \ b=\sum_{j \in B} d_j$, $ \ c=\sum_{j \in C} d_j$ and $d=\sum_{d \in D} |d_j|$. Then,
it is immediately seen that
\beqq
S(\lambda) &\geq & a \lambda^{\frac{n}{2}} \ + \ (b+c) \lambda^n -d \lambda^{\frac{9n}{10}}\\
&= &\left(\frac{d}{4} \lambda^{\frac{n}{2}} -d \lambda^{\frac{9n}{10}} + \frac{3d}{4} \lambda^n\right)
\ + \ \left( a-\frac{d}{4} \right) \left( \lambda^{\frac{n}{2}} - \lambda^n \right) \\
&& ~~~~~~~~~~~~\ + \ (b+c+a-d) \lambda^n .
\eeqq
To prove the assertion, it will be sufficient to prove the following three inequalities:
\be\label{b4}
\left \{\ba{l}
a+b+c\, -\, d \, > \,0,\\
4a>d,\\
\lambda^{\frac{n}{2}} \ + \ 3 \lambda^n -4 \lambda^{\frac{9n}{10}}\,\geq \,0, \quad \lambda \in (0,1].
\ea
\right .
\ee

Let us first begin to prove the first inequality in (\ref{b4}) which is obviously
equivalent to the inequality $\sum_{j=1}^md_j\,> \,0.$ We will prove this inequality by the following splitting
\beqq
\sum_{j=1}^m d_j &=& \sum_{j=1}^m \, \sum_{n-j+1 \leq r \leq m} f(r)\\
&=& \sum_{r=2}^m f(r) \, \sum_{n-r+1 \leq j \leq m} 1 \\
& = & \sum_{r=2}^m f(r) (r-1)\\
&=& \sum_{r=1}^mrf(r)\, -\, f(1)\, -\, \sum_{r=2}^m f(r)\\
&=& S_1 - S_2 - S_3.
\eeqq
Now,
$$S_1\, =\, \sum_{r=1}^m \left(r \sum_{l \leq \frac{n-1}{r}} \mu(l)\right)\, = \,\sum_{r=1}^m\sum_{rl \leq n-1} r \mu(l)
\, =\, \sum_{a=1}^m \sum_{rl=a} r \mu(l)
\,= \, \sum_{a=1}^m \varphi(a),
$$
where, as usual, $\varphi$ denotes Euler's totient function. From \cite{M} it is known that
\[\sum_{a=1}^m \varphi(a)\,\geq \,\frac{3m^2}{\pi^2}\,-\,\frac{1}{2}\,m\,\log\,m \,-\,\left(\frac{\gamma_0}{2}\,+\,\frac{5}{8}\right)m\,-\,1,
\]
where $\gamma_0$ denotes the Euler-Mascheroni constant. For the numbers $m$ under consideration, we may use the weaker estimate
\[\sum_{a=1}^m \varphi(a)\,\geq \,\frac{3m^2}{\pi^2}\,-\,\frac{1}{2}\,m\,\log\,m \,-\,m.
\]
Concerning $S_2$ and $S_3$, we have $S_2\,=\,M(m)\,\leq \,m,$ and
\[
S_3\,=\,\sum_{r=2}^mf(r)\,\leq \,\sum_{r=2}^m\frac{m}{r}\leq \, m\int_1^m\frac{d\,x}{x}\,=\,m\,\log \,m.\]
Hence
\[
\sum_{j=1}^md_j\,\geq \,\frac{3m^2}{\pi^2}\,-\,\frac{3}{2}\,m\,\log\,m \,-\,2m,
\]
which can easily be seen to be positive for the numbers $m$ in question.

%

Let us next prove the second inequality in (\ref{b4}), namely, $4a>d$.
Since
\[
|d_j|\,\leq \,\sum_{r=n-j+1}^m|f(r)|\,\leq \,\sum_{r=n-j+1}^m\frac{m}{r},
\]
we have for $j\in D$ the estimate
\begin{eqnarray*}
\sum_{j \in D} |d_j| &\leq &  \sum_{\frac{9n}{10} \leq j \leq n} \sum_{n-j+1 \leq r \leq n} \frac{m}{r}\nonumber\\
&=& m\sum_{r=1}^n \frac{1}{r} \sum_{\st{n-r+1 \leq j \leq n\\
0.9n \leq j \leq n}} 1\nonumber \\
&= &m \sum_{r \leq \frac{n}{10}} \frac{1}{r} \sum_{n-r+1 \leq j \leq n} 1 \ + \ m \sum_{0.1n \leq r \leq n} \frac{1}{r} \sum_{0.9n \leq j \leq n} 1\\
&\leq & m \sum_{r \leq {\frac{n}{10}}} 1 \ + \ m\sum_{\frac{n}{10} \leq r \leq n} \frac{1}{r} \frac{n}{10} \nonumber\\
&\leq & \frac{n^2}{10}\,+\,\frac{n^2}{10}\int_{\frac{n}{10}-1}^n\frac{d\,x}{x} \nonumber\\
&\leq &  \frac{n^2}{10} \left[\log\left(\frac{950}{85}\right) + 1\right ] \nonumber\\
&\leq & 0.342 n^2 <\,\frac{n^2-4n+3}{2}\,\leq \,4a .\nonumber
\end{eqnarray*}

Finally we prove the third inequality in (\ref{b4}), namely,
$$\lambda^{\frac{n}{2}} \ + \ 3 \lambda^n -4 \lambda^{\frac{9n}{10}} \,=\, \lambda^{\frac{n}{2}}
\left (1 \ + \ 3 \lambda^{\frac{n}{2}} -4 \lambda^{\frac{2n}{5}}\right ) \,\geq\,0 ~\mbox { for }~\lambda \in (0,1].
$$
 In order to prove this,
we let $x\,=\,\lambda^{\frac{n}{2}}$ and we see that it is sufficient to prove
\[
g(x)\,=\,1\,-\,4\,x^{\frac{8}{10}}\,+\,3x \,\geq\,0~\mbox { for }~ x\in (0,1].
\]
But this inequality is a direct consequence of the fact that $g(0)\,=\,1, g(1)\,=\,0$, and $g'(x)\,<\,0$ in $(0,1)$.
This completes the proof of the theorem.
\epf

\begin{rem}\label{rem1}
As a careful analysis of the proof reveals  we have actually proved a bit more, namely,
\[\sum_{j=0}^{n}\lambda^j\,-\,\sum_{k=1}^{n-1}\left\{M\left(\left[\frac{n-1}{k}\right]\right)
\left(\sum_{j=0}^{n-k}\lambda^j\right)\right\}\,\geq \,\lambda^n ~\mbox{ for }~\lambda \in [0,1].
\]
\end{rem}

Using again (\ref{f8}) we get another formulation of our theorem which is as follows.
\bthm\label{th2}
Let $n\in \IN\setminus\{1\}$ and $\lambda \in (0,1).$  Then
\[\sum_{k=1}^{n-1}\left(M\left(\left[\frac{n-1}{k}\right]\right)\lambda^{n-k+1}\right)\,>\,\lambda^{n+1}.\]
\ethm

\begin{rem}
In view of Remark \ref{rem1}, we can obtain a better inequality than Theorem \ref{th2}:
\[\sum_{k=1}^{n-1}\left(M\left(\left[\frac{n-1}{k}\right]\right)\lambda^{n-k+1}\right)\,>\,\lambda^{n}.\]
\end{rem}

\subsection*{Acknowledgement}
The authors thank the referee for his careful reading of the paper and his useful comments.

%

\end{document}